 \documentstyle[graphicx,amscd,amssymb,verbatim,12pt,righttag,color]{amsart}
 \newtheorem{theorem}{Theorem}[section]
 \newtheorem{Def}[theorem]{Definition}
 \newtheorem{Prop}[theorem]{Proposition}
 \newtheorem{Lem}[theorem]{Lemma}
 \newtheorem{Cor}[theorem]{Corollary}
 
 \newtheorem{Example}[theorem]{Example}

 \numberwithin{equation}{section}
 \renewcommand{\rm}{\normalshape}
 
\begin{document}

\title  {Spectral structure of digit sets \\ of self-similar tiles on ${\Bbb R}^1$}

\date{}
\author{Chun-Kit Lai}
\address{Department of Mathematics, The Chinese University of Hong Kong,  Hong Kong}
\author{Ka-Sing Lau}
\address{Department of Mathematics , The Chinese University of
Hong Kong , Hong Kong}
\author {Hui Rao}
\address {Department of Mathematics, Central China Normal University, Wuhan, China}
\email{cklai@@math.cuhk.edu.hk}
\email{kslau@@math.cuhk.edu.hk}
\email{hrao@@mail.ccnu.edu.cn}
\thanks {}

\date{March 31, 2011}
\thanks{ The research is supported in part by the HKRGC Grant and the Focused Investment Scheme of CUHK. }
\maketitle

\begin{abstract}
We study the structure of  the  digit sets ${\mathcal D}$ for the integral self-similar tiles $T(b,{\mathcal{D}})$ (we call such ${\mathcal D}$ a {\it tile digit set} with respect to $b$). So far the only available classes of such tile digit sets are the complete residue sets and the product-forms ([B], [K], [LW2], [LR]). Our investigation here is based on the spectrum of the mask polynomial $P_{\mathcal D}$, i.e., the zeros of $P_{\mathcal D}$ on the unit circle. By using the Fourier criteria of self-similar tiles of Kenyon [K] and Protasov [P], as well as the algebraic techniques of cyclotomic polynomial, we characterize the tile digit sets through some product of cyclotomic polynomials (kernel polynomials), which is a generalization of the product-form to higher order.
\end{abstract}

\bigskip

\begin{section} {\bf Introduction}
$\mbox{}$ Let $A$ be an $s\times s$  expanding matrix (i.e., all eigenvalues have moduli $> 1$) with integral entries and  $|\det A |= b$. Let also ${\mathcal{D}} \subset {\Bbb Z}^s$ with $\#{\mathcal D} =b$, we  call ${\mathcal{D}}$ a \textit{digit set}.  Consider the iterated function system (IFS) defined by
$\{S_j\}_{j=0}^{m-1}$ where
\begin{equation} \label {eq1.0}
S_j(x) = A^{-1}(x+d_j),    \qquad  0\leq j \leq m-1,
\end{equation}
 then these affine maps are contractive under certain norms in ${\Bbb R}^s$,  and it is well-known that there exists
a unique compact set $T := T(A,{\mathcal{D}})\subset {\Bbb R}^s$ satisfying the set-valued relation $AT = T+{\mathcal{D}}$.  Alternatively, $T(A,\mathcal{D})$ can be expressed as a set of radix expansions with base $A$ and digits in ${\mathcal{D}}$:
\begin{equation}\label{1.1}
T(A,{\mathcal{D}}) = \{\sum_{k=1}^{\infty}A^{-k}d_k : d_k\in
{\mathcal{D}}\}.
\end{equation}
 In [B], Bandt proved  that if  $|\det A| = b = \# {\mathcal D}$ and $T^o \neq \emptyset$, then $T(A, {\mathcal{D}})$ is a \textit{translational tile}. This means that there exists ${\mathcal{J}}\subset {\Bbb R}^s$ such that $\bigcup_{t\in{\mathcal{J}}}(T+t) ={\Bbb R}^s$, and the Lebesgue measure of $( T
+t) \cap (T +t')$ is zero for all  $t,t'\in{\mathcal{J}}$, $t\neq t'$. In that case we call  $T$ an {\it (integral) self-affine tile}, and a {\it self-similar tile} if the matrix $A$ is, in addition,  a scalar times an orthonormal matrix. The digit set ${\mathcal D}$ is called a {\it (self-affine) tile digit set} with respect to $A$, and ${\mathcal J}$ is called a \textit{tiling set} for $T$.

\medskip

 The study of self-affine tiles and the tiling theory was initiated by Thurston [Th] and Kenyon \cite{[K]}, and  the foundation and
basic properties were laid down in a series of papers by Lagarias and Wang [LW1-4]. By now there is a wealth of  literature  in this
subject,  with topics including the tiling theory, the geometric and
fractal structures,  the topological properties, the classification
problems and application to wavelet theory ([AT], [GH], [HLR],
[KL1,2], [LW1-4], [LL], [LR], [SW], [V]). The related topic of the
Fuglede's problem on  tiles and spectral sets has also received a
lot of attention recently ([Ko], [{\L}], [{\L}aW], [T]). However,
despite such intensive studies,  many of these aspects are still not
fully understood. Among them is the following basic question:

\medskip

{\bf (Q)}\  \ {\it Given a expanding integral matrix $A$, can we classify all the tile digit sets ${\mathcal{D}}$, and what is the structure of such ${\mathcal D}$?}

\medskip

The question turns out to be rather intriguing and challenging, it is not clear even in ${\Bbb R}^1$.  A well-known sufficient condition (it is also true in any dimension) for ${\mathcal D}$ to be a tile digit set is:\textit{ ${\mathcal{D}}$ is a complete residue set ({\rm mod} $A$)} \cite{[B]} (which is called {\it standard tile digit sets} [LW3]). This condition is also necessary in ${\Bbb R}^1$ if $b = p$ is a prime \cite{[K]}. It is also true in higher dimension when additional hypothesis is assumed ([LW3], \cite{[HL]}).

\medskip

The situation is far more complicated if $b$ is not a prime number. For example in ${\Bbb R}^1$ if $A = [b]$ with $b = 4$ and ${\mathcal{D}} = \{0,1,8,9\}$, then ${\mathcal D}$  is not a complete residue (mod $4$), but $T(4,{\mathcal{D}}) = [0,1]\cup[2,3]$ is a tile with a tiling set ${\mathcal J} = \{0,1\}\oplus 4{\Bbb Z}$.  In reference to the work of Odlyzko \cite{[O]},  Lagarias and Wang proposed a class of digit sets in ${\Bbb R}^s$ called \textit{product-form} \cite{[LW3]} as a generalization of the complete residue class. They also showed that for  $b = p^\ell$ a prime power,  then a tile digit set ${\mathcal{D}} $ of $b$  must be a ``product-form-like" set. More recently, Lau and Rao introduced a  \textit{weak product-form} and used it to classify all tile digit sets ${\mathcal{D}} (\subset {\Bbb Z})$ for $b = pq$, a product of two primes ([LR]).

\bigskip

In this paper, we only focus on self-similar tile in ${\Bbb R }^1$. We assume that $b>1$, $A = [b]$ and  $ \#{\mathcal D} =b$.  We  also assume, without loss of generality, that $0 \in {\mathcal D} \subset {\Bbb Z}^+$  and  g.c.d.$({\mathcal D})=1$. We
first recall that a cyclotomic polynomial $\Phi_d(x)$ is the
minimal polynomial of $e^{2\pi i /d}$. It is easy to see that  a product-form can be expressed into a product of cyclotomic polynomials (Section 3). Our
main result is to obtain a characterization of the tile digit sets as product of cyclotomic polynomials.   We call
$$
P_{{\mathcal D}}(x) = \sum_{d \in {\mathcal D}} x^d = 1+ x^{d_1} +\cdots + x^{d_{b-1}}
 $$
  a {\it mask polynomial} of ${\mathcal D}$. A well-known necessary and sufficient condition for $T(b, {\mathcal D})$ to be a self-similar tile, due to
Kenyon [K], is that the mask polynomial satisfies:  {\it for any
$m>0$, there exists $k$ such that
$$
P_{\mathcal D}(e^{2\pi i m/b^k})=0.
$$}
 The Kenyon criterion shows a close link of the zero set of the mask polynomials with the cyclotomic
polynomial.  Note that $P_{\mathcal D}(e^{2\pi i /d}) =0$ if and
only if $\Phi_d(x) |P_{\mathcal D}(x)$. We call the set \ $\{d>1:
\Phi_d (x) | P_{\mathcal D }(x)\}$ \ the {\it spectrum} of
${\mathcal D}$, and
$$
S_{\mathcal D} = \{p^\alpha >1:    p \ \hbox {prime},\
\Phi_{p^\alpha}(x) | P_{\mathcal D}(x) \}
$$
 the {\it prime-power spectrum}  of ${\mathcal D}$. They play  an important role in our consideration.

  \medskip

The spectrum was used extensively by Coven and Meyerowitz [CM] in their study of integer tiles $\mathcal A$ (i.e., $\#\mathcal A$ is finite, and  there exists $\mathcal L$ such that $\mathcal A \oplus {\mathcal L} = {\Bbb Z} $).  They gave an inductive characterization of such $\mathcal A$  when
 $\#{\mathcal A} =p^\alpha q^\beta, \ p, q$ are primes. Later,
{\L}aba and Wang [{\L}, {\L}aW] used the spectrum to investigate the  spectral sets and spectral measures. In the context of
self-similar tile digit sets $\mathcal D$, we known that ${\mathcal D}$ is an integer tile (see [LLR]). In here  we prove the following theorem  which characterizes the prime-power spectrum of ${\mathcal D}$ (Theorem
\ref{th4.2}).

\medskip

\begin{theorem} \label{th1.2}
 Let  $ b= p_1^{\alpha_1}\cdots p_k^{\alpha_k}$ be a product of prime powers and let ${\mathcal D}$ be a  tile digit set of $b$. Then
 $$
 S_{\mathcal D}= {\bigcup}_{j=1}^{\ k} S_{p_j}
 $$
where $S_{p_j} = \{p_j^{a_{j,1}}, \cdots , \ p_j^{a_{j, \alpha_j}}\}$ and $\{a_{j,1}, \cdots , a_{j, \alpha_j}\}$ is a complete residue set modulo $\alpha_i$.
\end{theorem}

\bigskip

It is well-known that if ${\mathcal D} = \{0, 1, \cdots, b-1\}$, then
 the mask polynomial
$$ P_{\mathcal D}(x) = 1 + x + \cdots +x^{b-1} = \prod_{d|b, d>1} \Phi_d(x);
$$
if ${\mathcal D}$ is a complete residue set modulo $b$, then
$
P_{\mathcal D} (x) = \Big (\prod_{d|b, d>1} \Phi_d(x) \Big )Q(x)
$
for some integral polynomial $Q(x)$. A direct generalization of this is

\bigskip

\begin{theorem} \label {th1.3} Suppose $\mathcal D$ satisfies

\vspace {0.2cm}

($ P_1$): \ \ {\it  for any factor $d (>1)$ of $b$, there exists $j(d) \geq 0$  such that $\Phi_d(x^{b^j}) | P_{\mathcal D}(x)$.}

\vspace {0.2cm}

\noindent Then ${\mathcal D}$ is a tile digit set and
\begin{equation} \label {eq1.2}
P_{\mathcal D} (x) = \Big (\prod_{d|b, d>1} \Phi_d(x^{b^{j(d)}})\Big ) Q(x)
\end{equation}
for some integral polynomial $Q(x)$.
\end{theorem}

It is easy to check that if  ${\mathcal D}$ is a product-form [LW3] or a weak product-form [LR], then ${\mathcal D}$ satisfies condition ($P_1$) (Corollary \ref{cor5.4}). Actually, condition ($P_1$) is satisfied for the more general {\it modulo product-form}, which is defined through the product form by allowing certain modulo action of each factor (see Definition \ref{Def5.5}). On the other hand, such simple expression of $P_{\mathcal D}$ in \eqref {eq1.2} is not sufficient to cover all tile digit sets.  For this we introduce the notion of $(P_k), k\geq 1$  and correspondingly, the more concrete \textit{higher order product-forms}. Roughly speaking,  a second order product-form can be constructed as follows:  we take a modulo product-form and rearrange the digits to form a product, then use this to construct a new product-form.  By doing so, we are able to construct new tile digit sets that are not known in literature.

\bigskip

By using a device of Protasov [P], we can represent the mask polynomial by a ``tree" of cyclotomic polynomials (a {\it $\Phi$-tree}) with respect to $b$, and a {\it blocking} ${\mathcal N}$ of the tree: a finite subset such that every path on the tree meets one and only one element of ${\mathcal N}$ (Section 5). Our main theorem is:

\medskip

\begin {theorem} \label {th1.4}
Let ${\mathcal{D}}$ be a set of non-negative integers with $\#{\mathcal D}=b$. Then the following are equivalent:

(i) ${\mathcal D}$ is a tile digit set;

(ii) there is a blocking ${\mathcal N}$ in the $\Phi$-tree such that
$$
P_{\mathcal D} (x) = \Big (\prod_{\Phi_d\in{\mathcal N}}\Phi_d(x)\Big ) Q(x);
$$

(iii) $P_{\mathcal{D}}(x)$ satisfies condition $(P_k)$  for some $k\geq 1$.

\end{theorem}

\bigskip

We call the product in (ii) a {\it kernel polynomial} of ${\mathcal
D}$. It plays a central role on the structure of the tile digit
sets.   For $b= pq$ with the prime-power spectrum equals
$\{p,q^{\ell}\}$, the digit set  ${\mathcal D}$ is a ($1^{st}$ order ) weak product-form ([LR]),  and the kernel polynomial is of the form $\Phi_p(x) \Phi_{q^\ell}(x^{p^{\ell}})$. For more
general case that $b=p^{\alpha}q$, we determine, in a separate paper [LLR], all the kernel polynomials that generate tile digit sets. By doing so, we  also show that all  tile digit sets are of modulo product-forms up to some order depending on the prime-power spectrum.

\medskip

For the organization of the paper, we first provide some basic
properties of the cyclotomic polynomials and integer tilings in Section 2; we also bring in the Kenyon criterion to study the mask polynomial $P_{\mathcal D}$ and prove Theorem \ref{th1.2}. In Section 3, We consider the condition ($P_1$) and the various kind of product-forms, and in Section 4, we set up the ($P_k$) condition and the higher order product-form.  In Section 5,
we make use of the Protasov's device to study the $\Phi$-tree
and prove Theorem \ref {th1.4}. We conclude in Section 6 with some remarks and open questions.

\end{section}

\bigskip
\bigskip

\begin{section}{\bf Preliminaries}

\bigskip

 In this section, we give a brief summary on the cyclotomic polynomials which is needed for later discussions.
We use $\Phi_d(x)$ to denote the \emph{$d$-th cyclotomic polynomial}. It is the minimal polynomial of the primitive $d$-th root of unity, i.e., $\Phi_d(e^{2\pi i/d}) = 0$. Note that the degree of $\Phi_d(x)$ is equal to the Euler-phi function $\varphi(d)$ (the number of relatively prime integers in $1,...,d-1$).  It is well known that
\begin {equation} \label {eq2.0}
x^n-1 = \prod_{d|n}\Phi_{d}(x).
\end{equation}
The formula provides a constructive way to find $\Phi_d$ inductively.  The following basic properties of cyclotomic polynomials will be used throughout the paper.

\bigskip

\begin{Prop}\label{Prop3.1}
The cyclotomic polynomials satisfies the following:

\vspace {0.2cm}

\ \  (i) If $p$ is a prime, then $\Phi_p(x) = 1+x+...+x^{p-1}$ and $\Phi_{p^{\alpha+1}}(x) = \Phi_p(x^{p^{\alpha}})$;

\vspace {0.2cm}

\ (ii) $\Phi_{s}(x^p) = \Phi_{sp}(x)$ if $p$ is prime and $p|s$, and

 \hspace {0.6cm} \ $\Phi_{s}(x^p) = \Phi_{s}(x)\Phi_{sp}(x)$ if $p$ is
prime but $p\nmid s$;

\vspace {0.2cm}

(iii) $\Phi_s(1) = \left\{
\begin{array}{ll}
0, & \hbox{if $s=1$}; \\
 p, & \hbox{if $s=p^{\alpha}$}; \\
  1, & \hbox{otherwise.}
 \end{array}
 \right.
 $

\end{Prop}

\bigskip

\begin{Prop}\label{Prop3.3} Let $b\geq 2$ be an integer. Then for any two distinct factors $d_1, d_2 \ (\neq 1)$ of $b$, and for any  integers  $u_1,u_2 \geq 0$, $\Phi_{d_1}(x^{b^{u_1}})$ and  $\Phi_{d_2}(x^{b^{u_2}})$ have no common factor.
\end{Prop}

\medskip
\noindent {\bf Proof.}\  Note that if $f(x)$ and $g(x)$ have no common divisor, then so do $f(x^n)$ and $g(x^n)$ (use $a(x)f(x)+b(x)g(x)=1$ for some polynomials $a(x),b(x)$).
Hence, we only need to show that $\Phi_{d_1}(x)$ and $\Phi_{d_2}(x^{b^{u_2-u_1}})$ (assuming $u_1\leq u_2$) have no common divisor. This is simple as $e^{2\pi i /d_1}$ is not a root of the polynomial  $\Phi_{d_2}(x^{b^k})$ for all $k\geq0$.
\qquad $\Box$

\bigskip

Let ${\Bbb Z}^+$ be the set of
non-negative integers. For ${\mathcal{A}}\subset {\Bbb Z}^+$,  we
let
  $$
  P_{\mathcal{A}}(x) = {\sum}_{a\in{\mathcal{A}}}x^a
  $$
   and call it the {\it mask polynomial} of ${\mathcal A}$.
  For a finite set ${\mathcal A} \subset {\Bbb Z}^+$, we use
\begin{equation}\label{3.2}
S_{\mathcal{A}} = \{ p^\alpha > 1 : \ p \ \hbox {prime},
 \ \Phi_{p^\alpha}(x)|P_{\mathcal{A}}(x)\}
\end{equation}
to denote the {\it prime-power spectrum} of ${\mathcal A}$, and $\widetilde S_{\mathcal A} = \{s>1: \ \Phi_s(x)|P_{\mathcal{A}}(x)\}$ the {\it spectrum } of ${\mathcal A}$.

We call a finite set ${\mathcal A} \subset {\Bbb Z}$ an {\it
integer tile} if there exists ${\mathcal L}$ such that ${\mathcal A}
\oplus {\mathcal L} = {\Bbb Z}$.    The class of integer tiles on ${\Bbb Z}$ has been studied in depth in connection with the factorization of cyclic groups and cyclotomic polynomials (\cite{[CM]}, \cite{[deB]}).   In [CM], Coven and Meyerowitz  introduced the following two conditions to study the integer tiles:

\vspace {0.2cm}
\noindent { (T1)} \  { \it $\#{\mathcal{A}} = P_{\mathcal{A}}(1)=\prod_{s\in
S_{\mathcal{A}}} \Phi_s(1)$,}

\vspace {0.1cm}

\noindent { (T2)} \ {\it For any distinct prime powers $s_1, \ldots, s_n \in S_{\mathcal A} $ , then
$s_1 \ldots s_n \in {\widetilde S}_{\mathcal{A}}$. }

\vspace {0.2cm}

 \noindent They showed that {\it if ${\mathcal{A}}\subset {\Bbb Z}^+$ is a finite set, and suppose $(T1)$ and $(T2)$ hold, then ${\mathcal{A}}$ tiles ${\Bbb Z}$ with period \ $n =l.c.m.(S_{\mathcal{A}})$.  Conversely, if $\mathcal{A}$ is an integer tile, then
$(T1)$ holds; if in addition  $\# {\mathcal A} = p^\alpha q^\beta, \ \alpha , \beta \geq 0$, then $(T2)$ holds.}
It is still an open question whether an integer tile must satisfy (T2) in general.
The class of integer tiles and the (self-affine) tile digit sets are closely related. In fact it is shown in [LLR] that \textit{all tile digit sets (in any dimension) are integer tiles}.

\bigskip

 In the following, we start to consider the tile digit sets of $b$ in ${\Bbb R}$. We assume that ${\mathcal{D}}\subset {\Bbb Z}^+$, \ $0\in{\mathcal{D}}$, and \ $g.c.d.({\mathcal{D}}) = 1$. For $T:= T(b, {\mathcal D})$,  we let $\chi_T$ denote the characteristic function of $T$, the self-similar identity $bT = T+{\mathcal{D}}$ can be expressed as a refinement equation
  $$
  \chi_T (x) = \sum_{d\in {\mathcal D}}\chi_T(bx - d), \mbox{ for a.e. $x\in{\Bbb R}$},
  $$
  and the Fourier transform is $\widehat \chi_T (\xi) = \prod_{k=1}^\infty P_{\mathcal D}(e^{2\pi i\xi /b^k})$. By using the Riemann-Lebesgue lemma, it is not difficult to derive from the product that for $\chi_T$ to be an $L^1$-function (i.e., $T$ is a self-similar tile since it has positive Lebesgue measure), then $\widehat \chi_T (m)=0$ for any integer $m \not =0$, so that one of the factors is $0$. The converse can be proved by considering certain tempered distribution (see e.g. [HL]).  This is the basic idea of the following criterion due to Kenyon (actually holds in ${\Bbb R}^s$).

\begin{theorem} {\rm [K]} \label{th4.1}\
 $T(b, {\mathcal D})$ is a  self-similar tile if and only if for each integer $m >0 $, there exists $k\geq 1$ (depending on $m$) such that
$$
P_{\mathcal{D}}(e^{2 \pi i {m}/{b^k}}) = 0.
$$
\end{theorem}

\bigskip

The factorization of $P_{\mathcal{D}}(x)$ is closely connected with the cyclotomic polynomials. By using the Kenyon criterion, we prove the following theorem which characterizes the prime power spectrum of ${\mathcal D}$.

\bigskip

\begin{theorem}\label{th4.2}
 Let  $b= p_1^{\alpha_1} \cdots p_k^{\alpha_k}$  be the product of prime powers and let $ {\mathcal D}$ be a  tile digit set of $b$. Then for each
$p_j$, and $0\leq \ell \leq \alpha_j-1$, there exists a unique $a_{j, \ell}\equiv \ell $ ({\rm {mod}} $\alpha_j$) such that
 $ \Phi_{p^{a_{j, \ell}}_j}(x)$ divides $\ P_{\mathcal{D}}(x)$, and the prime power spectrum of ${\mathcal D}$ is
 $$
 S_{\mathcal D} = {\bigcup}_{j=1}^k \{p_j^{a_{j,1}}, \cdots , \ p_j^{a_{j, \alpha_j}}\}.
 $$

\end{theorem}

\medskip

\noindent {\bf Proof}.
 We prove the theorem for $p_1$.
  For simplicity, we write $b = p_1^{\alpha_1} t$ \ and  $e(x) = e^{2\pi i x}$. For a fixed $0\leq \ell \leq \alpha_j-1$ and for each $n\geq 1$, we set
 $a_n = p_1^{\alpha_1-\ell} t^n.$
  By the Kenyon criterion, for any $n$, there exists $m =m(n)$ such that $e(a_n/{b^m})$
 is a root of $P_{\mathcal D}(x)$.

  \medskip

 We claim that there exists $n$ such that $n\geq m(n)$. Suppose otherwise, for any $n\geq 1$, we have $n<m(n)$ . Let  $u_n = p_1^{\ell+(m-1)\alpha_1}t^{m-n}$, then
  $ e(a_n/{b^m}) = e(1/u_n)$ are roots of $P_{\mathcal D}(x)$.  Note that $u_n > p_1^{\ell+(m-1)\alpha_1}$, hence $u_n$ tends to infinity as $n$ tends to infinity (since $n< m(n): =m)$. It means that $u_n$ contains infinitely many distinct numbers, it is a contradiction since $P_{{\mathcal D}}(x)$ is a finite degree polynomial.

  \medskip

  Hence, we can find $n$ such that $n\geq m(n)$ and
  $$
  e(a_n/{b^m}) = e(t^{n-m}/{p_1^{\ell+(m-1)\alpha_1}})
  $$
  is a root of $P_{\mathcal D}(x)$. Let $ a_{1, {\ell} } = \ell+(m-1)\alpha_1$, we have $\Phi_{p_1^{a_{1, \ell}}}(x)|P_{\mathcal{D}}(x)$. This shows the existence part in the theorem.

  \vspace {0.2cm}

 For the uniqueness, we note that we have  $ P_{\mathcal D}(1) = \#{\mathcal D} = b$. On the other hand, each $\Phi_{p_{j}^{a_{j,\ell}}}(1)=p_j$ (by Proposition \ref{Prop3.1}), and the product of all these $\Phi_{p_{j}^{a_{j,\ell}}}(1)$ is $  p_1^{\alpha_1}\cdots p_k^{\alpha_k} =b$. Moreover, $\Phi_{p_{j}^{a_{j,\ell}}}(x)$ divides $P_{\mathcal D}(x)$. It follows that the set of all these ${p_{j}^{a_{j,\ell}}}$ is precisely the prime power spectrum $S_{\mathcal D}$ and the uniqueness follows.
 \qquad $\Box$

 \bigskip

  It follows from \eqref {eq2.0} that  if  ${\mathcal{D}}$ is a  complete residue set modulo $b$ ({\it standard digit set} [LW3]), then
$
P_{\mathcal D} (x) = \Big (\prod_{d|b, d>1} \Phi_d(x) \Big )Q(x)
$
for some integral polynomial $Q(x)$.  We see that the $S_{\mathcal{D}}$ in Theorem \ref{th4.2} is simply the set $p_j^i$, with $i=1,\cdots,\alpha_i$. The theorem implies that if we change to the other tile digit sets of $b$, the prime power factors are preserved in a modulo way.

\bigskip

\begin{Example}
{\rm Let ${\mathcal D}=\{0,1,8,9\}$ and $b=4=2^2$. It is known that $T(b,{\mathcal D}) =[0,1]\cup[2,3]$ is a tile. Its mask polynomial is
$$
P_{\mathcal D}(x) = 1+x+x^8+x^9 = \Phi_2(x)\Phi_{2^4}(x).
$$
We see that $\{1,4\}$ is a complete residue (mod 2).

\medskip

 If ${\mathcal{A}} = \{0,1,4,5\}$ (it is an integer tile), we have
 $$
 P_{\mathcal{A}}(x) = 1+x+x^4+x^5 = \Phi_{2}(x)\Phi_{2^3}(x).
 $$
 Note that $\{1,3\}$ is not a complete residue (mod 2), this shows that ${\mathcal A}$ is not a tile digit sets by Theorem \ref{th4.2}.}
\end{Example}
\end{section}

\bigskip
\bigskip

 \begin {section} {\bf Modulo product-forms}

 \medskip

    It is known that if $b$ is a prime, then $\mathcal D$ is a tile digit set if and only if  ${\mathcal D} \equiv {\Bbb Z}_b$ \cite{[K]}. But when $b$ is not a prime, the problem is vastly more complicated. The following classes of ${\mathcal D}$ are some of the basic tile digit sets that are known.

\medskip

 \begin {Def} \label{Def5.1} ${\mathcal{D}}$ is
called a product-form digit set (with respect to $b$) \cite{[LW3]} if
\begin{equation} \label {eq5.0}
{\mathcal{D}} = {\mathcal{E}}_0 \oplus b^{l_1} {\mathcal{E}}_1 \oplus \ldots \oplus b^{l_k}
{\mathcal{E}}_k
\end{equation}
where
$  {\mathcal{E}} ={\mathcal{E}}_0 \oplus{\mathcal{E}}_1 \oplus \ldots \oplus{\mathcal{E}}_k \equiv {\Bbb{Z}}_b $, and $0\leq l_1 \leq
l_2\leq \ldots \leq l_k$;  if ${\mathcal{E}} = \{0,1,2,\ldots ,b-1\}$, then ${\mathcal{D}}$ is called a strict product-form  \cite{[O]}.

\vspace{0.2cm}
Furthermore if ${\mathcal D}'\equiv {\mathcal D}$ ({\rm mod}   $b^{l_k+1})$ for ${\mathcal D}$ in (\ref {eq5.0}),  then ${\mathcal D}'$ is called a weak product-form \cite{[LR]}.

\end{Def}

\bigskip

In terms of the mask  polynomials,  a product-form can be expressed as
\begin{equation}\label {5.1}
P_{\mathcal{D}}(x)= P_{{\mathcal{E}}_0}(x)P_{{\mathcal{E}}_1}(x^{b^{l_1}})\ldots
P_{{\mathcal{E}}_k}(x^{b^{l_k}}),
\end{equation}
where
\begin{equation} \label {5.2}
P_{\mathcal{E}}(x) = P_{{\mathcal{E}}_0}(x)\ldots P_{{\mathcal{E}}_k}(x)
\equiv 1+x+\ldots +x^{b-1}({\mbox{mod}}\ x^b-1).
\end{equation}
For the weak product-form,
\begin{equation} \label {5.3}
P_{\mathcal{D}}(x)= P_{{\mathcal{E}}_0}(x)P_{{\mathcal{E}}_1}(x^{b^{l_1}})\ldots
P_{{\mathcal{E}}_k}(x^{b^{l_k}})+ (x^{b^{l_k+1}} -1) Q(x)
\end{equation}
for some integral polynomial $Q(x)$.

\bigskip

In view of the  expressions in \eqref {5.1}-\eqref {5.3},  we introduce a new condition on $\mathcal D$:

 \medskip

 \noindent {$({\bf P}_1)$} \  {\it  for any $d>1$ and \ $d|b$, there exists $j\geq 0$ (depends on $d$) such that $\Phi_d(x^{b^j}) | P_{\mathcal D}(x)$.}

 \vspace {0.2cm}
\noindent We also say that a polynomial $P(x)$ satisfies $({P}_1)$ if the above holds with  $P_{\mathcal D}(x)$ replaced by  $P(x)$.

\bigskip

\begin{theorem} \label{Prop5.3}
Suppose ${\mathcal{D}}$ satisfies condition $(P_1)$, then ${\mathcal{D}}$ is a tile digit set, and
\begin {equation}\label {5.5}
P_{\mathcal D}(x) = \prod_{d|b, d>1} \Phi_d(x^{b^{j(d)}}) \ Q(x)
\end{equation}
where $j(d)$ is an integer depending on $d$.
\end{theorem}

\medskip

\noindent {\bf Proof}. \
 Let $m>0$,  we claim that  there exists $d|b,\  d>1$ such that $\Phi_d(e^{2\pi i m/b^k})=0$. Consider the standard tile digit set
${\mathcal{E}}= \{0,1,\ldots,b-1\}$, by Theorem \ref{th4.1},  there exists $k\geq1$ such that
$$
 P_{\mathcal{E}}(e^{2\pi i {m}/{b^k}})  = 0.
 $$
Let $a= g.c.d.(m,b^k)$ and let $s = m/a,\  \ell= b^k/a$, then $s, \ \ell \ (>1)$ are relatively primes and $P_{\mathcal{E}}(e^{2\pi i {s}/{\ell}})  = 0$. It follows that
\begin {equation} \label {5.4}
\Phi_{\ell}(x)|P_{\mathcal{E}}(x).
\end{equation}
As $P_{\mathcal E}(x) = \prod_{d|b,  d >1}\Phi_{d}(x)$,
we have $\ell = d$ for some $d>1$ and the claim follows.

\vspace{0.2cm}

By condition $(P_1)$, there exists $j:=j(d) \geq 0$ such that
\begin{equation} \label {5.6}
P_{\mathcal{D}}(x) = \Phi_{d}(x^{b^j})Q_d(x)
\end{equation}
for some polynomial $Q_d(x)$. We hence have
$$
P_{\mathcal{D}}(e^{2\pi i {m}/{b^{k+j}}})\
 = \  \Phi_d(e^{2\pi i {m}/{b^k}})Q_d(e^{2\pi i {m}/{b^{k+j}}}) \
= \ 0.
$$
 This implies ${\mathcal{D}}$ satisfies the Kenyon criterion and hence it is a tile digit set.

 For the expression in \eqref{5.5}, we only need to use \eqref {5.6}  and  observe that those  $\Phi_d(x^{b^{j(d)}})$'s are distinct factors of $P_{\mathcal D}(x)$ (by Proposition \ref {Prop3.3}). \qquad $\Box$

 \bigskip

 \begin {Cor} \label{cor5.4} A product-form (or a weak product-form) digit set is a tile digit set.
 \end{Cor}

 \medskip

\noindent  {\bf Proof}.  It follows from \eqref {5.2} that
$$
P_{{\mathcal{E}}_0}(x)\ldots P_{{\mathcal{E}}_k}(x) = P_{\mathcal{E}}(x) =  1+x + \ldots +x^{b-1} + (x^b-1)Q(x)
$$
for some polynomial $Q(x)$. This implies that for any $d|b$ and $d>1$, $\Phi_d(x)|P_{{\mathcal{E}}_i}(x)$ for some $i=0,\dots,k$. Hence $\Phi_d(x^{b^{l_i}})|P_{{\mathcal{E}}_i}(x^{b^{l_i}})$. On the other hand, by \eqref{5.1},
$$
P_{\mathcal{D}}(x)= P_{{\mathcal{E}}_0}(x)P_{{\mathcal{E}}_1}(x^{b^{l_1}})\ldots P_{{\mathcal{E}}_k}(x^{b^{l_k}}),
$$
we obtain $\Phi_d(x^{b^{l_i}})|P_{\mathcal{D}}(x)$, and hence ${\mathcal D}$ is a tile digit set by Proposition  \ref{Prop5.3}.

\vspace{0.2cm}

For the weak product-form, we need to observe that the last factor in  \eqref{5.3} is divisible by $\Phi_d(x^{b^{l_i}})$.
 \qquad $\Box$

\bigskip

We now introduce a more general kind of digit sets that satisfies $(P_1)$.  Consider the product-form ${\mathcal D}$ in Definition \ref{Def5.1}, we define a decomposition of the spectrum: let $S_i =\{d >1: d|b, \ \Phi_d(x)|P_{{\mathcal{E}}_i}(x)\}$ and let
\begin{equation}\label{5.7--}
\Psi_i(x) = \prod_{d\in S_i}\Phi_{d}(x).
\end{equation}
Then  $\Psi_i(x)|P_{{\mathcal{E}}_i}(x)$, hence $\Psi_i(x^{b^{l_i}})|P_{{\mathcal{E}}_i}(x^{b^{l_i}})$. Let
\begin{equation}\label{5.7-}
K^{(i)}_{{\mathcal {D}}}(x) = \Psi_0(x)\Psi_1(x^{b^{l_1}})...\Psi_i(x^{b^{l_i}}), \quad 0 \leq i\leq k \ ,
\end{equation}
and denote $ K^{(k)}_{{\mathcal {D}}}(x)$ by $K_{\mathcal D}(x)$. Then $K_{{\mathcal D}}(x)|P_{{\mathcal D}}(x)$. Moreover, $K_{{\mathcal D}}(x)$ satisfies $(P_1)$. It is because for any $d>1$ a factor of $b$, we have $d\in{\mathcal S}_i$ for some $i$, hence $\Phi_{d}(x^{b^{l_i}})|\Psi_i(x^{b^{l_i}})$ so that $\Phi_d(x^{b^{l_i}})|K_{{\mathcal D}}(x)$ also.

\medskip

\begin{Def}\label{Def5.5} We say that ${\mathcal D}$ is a modulo product-form if  ${\mathcal D} := {\mathcal D}^{(k)}$ is defined through a product-form ${\mathcal{D}}'  = {\mathcal{E}}_0 \oplus b^{l_1} {\mathcal{E}}_1 \oplus \ldots  \oplus b^{l_k}
{\mathcal{E}}_k,  \ 0\leq l_1  \ldots \leq l_k $ as follows: let
$$
n_i  \ = \ l.c.m.\ \{s: \ \Phi_{s}(x)\ |\ K_{{\mathcal D}'}^{(i)}(x)\},
$$
and
\begin{equation} \label{5.8}
\left\{
  \begin{array}{ll}
    {\mathcal{D}}^{(0)} \equiv {\mathcal{E}}_0 \ ({\mbox{mod}} \ n_0),\\
    \vdots \\
    {\mathcal{D}}^{(i)} \equiv {\mathcal{D}}^{(i-1)} \oplus b^{l_i}{\mathcal{E}}_i \ ({\mbox{mod}} \ n_i),\\
    \vdots\\
    {\mathcal{D}}^{(k)} \equiv {\mathcal{D}}^{(k-1)} \oplus b^{l_k}{\mathcal{E}}_k \ ({\mbox{mod}} \ n_k)
  \end{array}
\right.
\end{equation}

\end{Def}

\bigskip

\noindent {\bf Remarks: }

(1) Note that $b^{l_i}|n_i$ and  $n_i | b^{l_i+1}$. Indeed  let ${\mathcal R}_i = \{s: \Phi_{s}(x)|K^{(i)}_{{\mathcal D}'}(x)\}$, and  let $d \in S_i$, then $\Phi_{d}(x^{b^{l_i}})|\Psi_{i}(x^{b^{l_i}})$. Using Proposition \ref{Prop3.1}(ii), we have
$ db^{l_i}\in {\mathcal R}_i$.
This means  $db^{l_i}|n_i$, therefore $b^{l_i}|n_i$.

Also for any $s \in {\mathcal R}_i$, $\Phi_{s}(x)|\Phi_{d}(x^{b^{l_j}})$ for some $0\leq j\leq i$ and $d \in S_j$. Since $d |b$, we have $\Phi_{d}(x^{b^{l_j}}) |  (x^{b^{l_i+1}}-1)$. Hence $\Phi_{s}(x)| (x^{b^{l_i+1}}-1)$, so that $s|b^{l_i+1}$. The definition of \textit{l.c.m.} implies $n_i|{b^{l_i+1}}$.

\medskip

(2)  It is clear that a product-form is a modulo product-form by ignoring all the modulo actions, and so is the weak product-form by keeping the last modulo.

\medskip

\bigskip

\begin{theorem}\label{th5.6}
Let ${\mathcal{D}}$ be a modulo product-form derived from the product form ${\mathcal D}'$,  then the polynomial $K_{{\mathcal D}'}(x)$ in \eqref{5.7-} divides  $P_{{\mathcal D}}(x)$,  hence ${\mathcal{D}}$ satisfies condition $(P_1)$ and is a tile digit set.
\end{theorem}

\bigskip

\noindent {\bf Proof}.  For each $i =0,...,k$, by the definition of $n_i$, $x^{n_i}-1 = \prod_{d|n_i}\Phi_d(x)$ and $\Psi_{i}(x^{b^{l_i}})$ is a product of cyclotomic polynomial, we have $\Psi_{i}(x^{b^{l_i}})\ |\ (x^{n_i}-1)$. Moreover, we know $\Psi_{i}(x^{b^{l_i}})|P_{{\mathcal E}_i}(x^{b^{l_i}})$. Hence in view of
$$
P_{{\mathcal{D}}^{(i+1)}}(x) =P_{{\mathcal{D}}^{(i)}}(x)P_{{\mathcal{E}}_{i}}(x^{b^{l_i}})+(x^{n_i}-1)
Q_{i+1}(x),
$$
we conclude that $\Psi_{d}(x^{b^{l_i}})\ |\ P_{{\mathcal{D}}^{(i+1)}}(x)$.

Note that the definition of $n_{i+1}$ implies that $\Psi_{i}(x^{b^{l_i}})\ |\ (x^{n_{i+1}}-1)$; this together with the expression of $P_{{\mathcal{D}}^{(i+2)}}(x) $ implies that $\Psi_{i}(x^{b^{l_i}})\ |\ P_{{\mathcal{D}}^{(i+2)}}(x)$. Continuing the process, we obtain $\Psi_{i}(x^{b^{l_i}})\ |\ P_{\mathcal{D}}^{(k)}(x)$. This shows that  $K_{{\mathcal D}'} = K^{(k)}_{{\mathcal D}'}$ in \eqref{5.7-} divides  $P_{\mathcal D}$. As $K_{{\mathcal D}'}$ satisfies $(P_1)$, we have ${\mathcal D} = {\mathcal D}^{(k)}$ satisfies $(P_1)$ also, and ${\mathcal D}$ is a tile digit set following from Theorem \ref{Prop5.3}.  \qquad $\Box$

\bigskip

We end this section by giving an example to illustrate the construction of a modulo product-form digit set.

\bigskip

\begin {Example}\label{ex5.6}   {\rm Let $b=12$ and write ${\mathcal E} = {\mathcal E}_0\oplus {\mathcal E}_1 \oplus {\mathcal E}_2 = \{0,1\} \oplus \{0,4, 8\}\oplus \{0,2\}$. Let
$$
{\mathcal D}' = {\mathcal E}_0\oplus {\mathcal E}_1 \oplus 12 {\mathcal E}_2
$$
Then ${\mathcal D'}$ is a (strict) product-form digit set. Let ${\mathcal D}^{(0)} = \{0,1\}$, then $n_0=2$.  To construct ${\mathcal D}^{(1)}$, we first evaluate $n_1$. Since $1+x^4+x^8=\Phi_3(x^4)$, we have
$$
K_{{\mathcal D}'}^{(1)}(x) = \Psi_0(x) \Psi_1(x) = \Phi_2(x)\Phi_3(x^4)= \Phi_2(x)\Phi_3(x)\Phi_6(x)\Phi_{12}(x).
$$
Hence  $n_1 =12$. We can choose
$$
{\mathcal D}^{(1)} = \{0,1, 4, 8, 9, 17 \} \ \equiv \ {\mathcal D}^{(0)} \oplus {\mathcal E}_1  \ \ (\hbox {mod}\  12).
$$
Next we observe that $K_{{\mathcal D}'}^{(2)}(x)= \Psi_0(x) \Psi_1(x)\Psi_2(x^{12})$ is given by
$$
\Phi_2(x)\Phi_3(x^4)\Phi_{2^2}(x^{12})=
\Phi_2(x)\Phi_3(x)\Phi_6(x)\Phi_{12}(x) \Phi_{16}(x) \Phi_{48}(x).
$$
This implies $n_2 =48$. We let
$$
{\mathcal D} = \{ 0,1,4,8,9,17,25,33,41,72,76,80 \} \equiv \ {\mathcal D}^{(1)} \oplus 12{\mathcal E}_2  \ \ (\hbox {mod} \ 48).
$$
Then ${\mathcal D}$ is a modulo product-form.

\medskip

On the other hand, ${\mathcal D}$ is not a product-form.
 For if so, then ${\mathcal{D}}$  can only be of the form ${\mathcal{E}}_0\oplus12{\mathcal{E}}_1$ with ${\mathcal{E}}_0\oplus {\mathcal{E}}_1 \equiv {\Bbb Z}_{12}$.  It is direct to check that
${\mathcal{E}}_0 = \{0, 1, 4, 8, 9, 17\}$ necessarily. But then it cannot be written in the product-form as needed.

\vspace {0.2cm}

It is clear that  ${\mathcal D}$ is not a weak product-form either, since taking modulo $12^n$ with $n\geq2$ reduces back to the original set which is not a product-form, while taking modulo $12$ contains only $6$ digits.}\end{Example}

\end{section}

\bigskip
\bigskip

\begin{section}{\bf Higher order product-forms}

The condition $(P_1)$ and the modulo product-form in the last section does not cover all tiles digit sets (see Example \ref {ex6.2}).  In this regard, we set up a higher order analog in this section, which will be studied in detail in the following sections.

\medskip

Let $\Phi_d(x^{b^{j}})$ be as in the definition of $(P_1)$, then by Proposition \ref {Prop3.1}(ii),
$$
\Phi_d(x^{b^{j}}) = \Phi_{t_1}(x) \ldots \Phi_{t_n}(x),
$$
and  each $\Phi_{t_j}(x)$ is a factor of $P_{\mathcal D}(x)$ (by Proposition \ref{Prop5.3}).  We observe that  a more relaxed condition that $\Phi_{t_j}(x^{b^{i_j}})$ is a factor of $ P_{\mathcal D}(x)$ also suffices for the Kenyon criterion to hold. We formulate this in the following. For $d|b$, let
\begin{equation} \label {eq5.2}
j_1 = j_1(d) :=
\min \{j: \exists  \mbox { a factor } \Phi_t(x) \ \mbox {of } \Phi_{d}(x^{b^j})  \backepsilon  \Phi_t(x)|P_{\mathcal{D}}(x) \}
\end{equation}
(and $ j_1= \infty$ if no factor $\Phi_t(x)$ exists). Define
\bigskip

(${\bf P}_2$)  \ \ {\it  For each $d|b, d>1$,  $j_1(d)<\infty $ and for any factor $\Phi_{t_1}(x)$ of $\Phi_d(x^{b^{j_1}})$, there exists $j_2 \geq 0$ (depends on $t_1$) with $\Phi_{t_1}(x^{b^{j_2}})\  |\  P_{\mathcal D}(x)$}.

\bigskip

Likewise we can repeat the same procedure of defining $j_2 = j_2(t_1)$ and  find $j_3$ for the factors of $\Phi_{t_1}(x^{b^{j_2}})$ to define $(P_3)$, and inductively for  $(P_k), k \geq 3$.

\bigskip

It is clear that $(P_1)\Rightarrow (P_2)$ by putting $j_1 =j$  and $j_2 = 0$. Also $(P_{k-1})\Rightarrow (P_k)$.

\bigskip

\begin{Prop}\label{Prop6.1}
Suppose ${\mathcal{D}}$  satisfies $(P_k)$, then ${\mathcal{D}}$ is a  tile digit set.
\end{Prop}

\medskip

\noindent {\bf Proof}.  We will check for  the case $(P_2)$ that the Kenyon criterion is fulfilled, the general case follows from the same idea. As in Theorem \ref {Prop5.3}, for $m>0$, there exists   ${\ell}_1>1$  with \ ${\ell}_1|b$ \ and \  $\Phi_{{\ell}_1} (e^{2 \pi i m/b^k}) =0$.
Take $j_1$ as in the assumption of $(P_2)$,  note that
$$
 \Phi_{{\ell}_1}((e^{2\pi i {m}/{b^{j_1+k}}})^{b^{j_1}}) =  \Phi_{{\ell}_1}(e^{2\pi i {m}/{b^{k}}}) = 0.
$$
Let $a = {g.c.d.(m,b^{j_1+k})}$ and let ${\ell}_2 = {b^{j_1+k}}/a, \ c = m/a$; they are relatively prime. From the above, we have  $\Phi_{{\ell}_1}((e^{2 \pi i c/l_2})^{b^{j_1}}) =0$, which implies that  $\Phi_{{\ell}_2}(x)|\Phi_{l_1}(x^{b^{j_1}})$. The $(P_2)$ assumption implies that there exists $j_2\geq 0$ such that $P_{{\mathcal{D}}}(x) = \Phi_{{\ell}_2}(x^{b^{j_2}})Q(x)$ for some polynomial $Q(x)$. Let  $j_1+j_2+k\geq 1$, then we have
$$
P_{{\mathcal{D}}}(e^{2\pi i {m}/{b^{j_1+j_2+k}}}) = \Phi_{{\ell}_2}(e^{2\pi i m /{b^{j_1+k}}})Q(e^{2\pi i {m}/{b^{j_1+j_2+k}}}) =0.
$$
This verifies the Kenyon criterion for ${\mathcal{D}}$,  and hence ${\mathcal{D}}$ is a  tile digit set.
\qquad $\Box$

\bigskip

 Next we give a concrete class of digit sets that satisfies  the $(P_k)$ condition. We regard the product-form and modulo product-form in  Definition \ref{Def5.1}  and \ref{Def5.5} as $1^{\rm st}$-order,  we define

\bigskip

 \begin{Def}\label{Def6.3}
${\mathcal{D}}$ is
called a $2^{\rm nd}$-order product-form (with respect to $b$) if
\begin{equation}\label{eq5.3}
{\mathcal{D}} = {\mathcal G}_0 \oplus b^{l_1} {\mathcal G}_1 \oplus \ldots \oplus b^{l_k} {\mathcal G}_k
\end{equation}
where $0\leq l_1 \leq
l_2\leq \ldots \leq l_k$, ${\mathcal G}={\mathcal G}_0 \oplus{\mathcal G}_1 \oplus \ldots \oplus{\mathcal G}_k$, and
${\mathcal G} $ itself is a modulo product-form as in Definition \ref{Def5.5} (possibly in another decomposition different from the ${\mathcal G}_i$).

For the above ${\mathcal G}$,  let $ S^{(2)}_i = \{s: \Phi_s(x)|P_{{\mathcal G}_i}(x), \ \Phi_s(x)|K_{{\mathcal G}}(x) \}$ where  $K_{{\mathcal G}}$ is defined in \eqref{5.7-},  we use the same procedure as in Definition \ref{Def5.5} to define the $2^{\rm nd}$-order modulo product-form and the corresponding $K_{{\mathcal D}}$. Inductively we can define the  $k^{\rm th}$-order modulo product-form.
\end{Def}

\medskip

\noindent {\bf Remark}. \ Roughly speaking, we can produce new digit sets as follows: we start with a modulo product-from ${\mathcal G}$, rearrange its digits to form a product, then use it to construct the  $2^{ \rm nd}$-order product-form, and then the $2^{\rm nd}$-order modulo product-form.

\bigskip

\begin{theorem}\label{th6.4}
Let ${\mathcal{D}}$ be a $k^{\rm th}$-order product-form (or a $k^{\rm th}$-order modulo product-form),  then it satisfies the $(P_k)$ condition and hence a tile digit set.
\end{theorem}

\noindent {\bf Proof}.
We only prove the case when ${\mathcal{D}}$ is a $2^{\rm nd}$-order product-form,  the other cases are similar. Let  $d|b$ and $d>1$, since ${\mathcal G}$ satisfies $(P)$, by Theorem \ref{Prop5.3}, there exists $j$ such that $\Phi_{d}(x^{b^{j}})|P_{\mathcal G}(x)$. As $P_{{\mathcal G}}(x) = P_{{\mathcal G}_0}(x)\ldots P_{{\mathcal G}_k}(x)$, we let $i$ be the first index such that there exists a factor $\Phi_t(x)$ of $\Phi_{d}(x^{b^{j}})$ divides $P_{{\mathcal G}_i}(x)$. Let $j_1$ be as in \eqref {eq5.2} corresponding to ${\mathcal D}$, it follows from \eqref{eq5.3} that $j_1 = j_1(d) = j+l_i$.

\medskip

Now for each  factor $\Phi_e(x)$ of $\Phi_{d}(x^{b^{j_1}})$, by
 $$
 \Phi_{d}(x^{b^{j_1}}) = \prod_{\Phi_{e'}(x)|\Phi_{d}(x^{b^j})}\Phi_{e'}(x^{b^{l_i}})\ ,
 $$
 we have $\Phi_e(x)|\Phi_{e'}(x^{b^{l_i}})$ for some $e'$ in the above product.  In view of $\Phi_{d}(x^{b^{j}})|P_{\mathcal G}(x)$, there exists $m$ such that $\Phi_{e'}(x)|P_{{\mathcal G}_{m}}(x)$. By the choice of $i$, we have $m\geq i$. Let $j_2 = j_2(e) = l_{m}-l_i$, then
 $$
 \Phi_e(x^{b^{j_2}})\ |\ \Phi_{e'}(x^{b^{l_{m}}}) \quad \mbox{and} \quad \Phi_{e'}(x^{b^{l_{m}}})\ |\ P_{{\mathcal G}_{m}}(x^{b^{l_{m}}}).
 $$
 As $P_{{\mathcal G}_{m}}(x^{b^{l_{m}}})|P_{{\mathcal{D}}}(x)$, this shows that $\Phi_e(x^{b^{j_2}})|P_{{\mathcal{D}}}(x)$.
\qquad $\Box$

\bigskip

The following diagram indicates the implications of the new classes of tile digit sets for a given $b$.
$$
\begin{array}{ccccccc}
  1^{\rm st}\hbox {-order  mpf }& \Rightarrow & 2^{\rm nd}\hbox {-order  mpf } & \Rightarrow & \ldots &  &   \\
  \Downarrow &  & \Downarrow &  &  &  &   \\
  (P_1) & \Rightarrow & (P_2) & \Rightarrow & \ldots & \Rightarrow  &{\mathcal{D}} \hbox { is  a  tile digit set}.
\end{array}
$$
(mpf means modulo product-form.) The $(P_k)$ can be understood better through the tree structure of cyclotomic polynomials to be developed in the next section;  we show that the converse of the last implication also holds, i.e., every tile digit set of $b$ must satisfy condition $(P_k)$ for some $k$.  On the other hand, it is still unclear if a digit set satisfies  $(P_k)$ for some $k$, then the digit set must be $k^{\rm th}$-order product-form.

\bigskip

We conclude this section by constructing a non-trivial $2^{\rm nd}$-order product-form digit set.

\begin{Example}\label{ex6.2}
 Let $b =12$ and
$$
{\mathcal{D}} = \{0,1\}\oplus2^4\{0,6\}\oplus 2^7\cdot3^2\{0,2,4\}
$$
Then ${\mathcal D}$ is a $2^{\rm nd}$-order product-form and it satisfies $(P_2)$, but not $(P_1)$.
\end{Example}

Indeed by rearranging the terms, we can write  ${\mathcal D}$  as
$$
{\mathcal{D}} = \{0,1\}\oplus12\{0,8\}\oplus (12)^2\{0,16,32\}: = {\mathcal G}_0\oplus12{\mathcal G}_1\oplus(12)^2{\mathcal G}_2.
$$
It is clear that  ${\mathcal G}: = {\mathcal G}_0\oplus{\mathcal G}_1\oplus{\mathcal G}_2 = \{0,1,8,9,16,17\}\oplus12\{0,2\}$ is a $1^{st}$-order product-form. Therefore ${\mathcal D}$ is a $2^{\rm nd}$-order product-form.

\vspace {0.2cm}
On the other hand, $\mathcal D$ does not satisfy $(P_1)$. Suppose otherwise, consider
\begin{equation}\label{6.1}
\begin{aligned}
P_{\mathcal{D}}(x) =& \Phi_2(x)\Phi_{2}(x^{6\cdot2^4})\Phi_3(x^{2^8\cdot3^2})\\
=& \Phi_2(x)\cdot \Phi_{2^6}(x)\Phi_{2^6\cdot3}(x) \cdot \Phi_{3^3}(x)
\Phi_{2\cdot3^3}(x) \ldots \Phi_{2^8\cdot3^3}(x)\\
\end{aligned}
\end{equation}
and for the factor $d=4$, there exists $j\geq 0$ such that $\Phi_4(x^{12^j})|P_{\mathcal{D}}(x)$. Since $\Phi_{2^{2j+2}}(x)|\Phi_4(x^{12^j})$, hence $\Phi_{2^{2j+2}}(x)|P_{\mathcal{D}}(x)$. From (\ref{6.1}), we must have $2j+2 = 6$, so that $j=2$. Thus $\Phi_4(x^{12^2})|P_{\mathcal{D}}(x)$. But this is impossible because
$$
\Phi_4(x^{12^2}) = \Phi_{2^6}(x)\Phi_{2^6\cdot3}(x)\Phi_{2^6\cdot3^2}(x),
$$
and  $\Phi_{2^6\cdot3^2}(x)$ does not divide $P_{\mathcal{D}}(x)$. \qquad $\Box$


\medskip
 It is also interesting to observe the following:

\begin{Example}\label{ex6.3} If we multiply one more factor  $\Phi_{2^6\cdot3^2}(x)$ to (\ref{6.1}), then it is the mask polynomial for another tile digit set ${\mathcal D} = \{0,1\} \oplus 12^2\{0, 2\} \oplus 12^2\{0, 16, 32\}$, which is a $1^{\rm st}$-order product-form of $b=12$.
\end{Example}

The expression of $\mathcal D$ follows from
\begin {eqnarray*}
P_{{\mathcal{D}}}(x) & = &\Phi_2(x)\cdot \Phi_{2^6}(x)\Phi_{2^6\cdot3}(x)\Phi_{2^6\cdot3^2}(x) \cdot \Phi_{3^3}(x)
\Phi_{2\cdot3^3}(x) \ldots \Phi_{2^8\cdot3^3}(x)\\
&=& \Phi_{2}(x)\Phi_{2}(x^{12^2\cdot 2})\Phi_{3}(x^{12^2\cdot 2^4}) \ .
\end {eqnarray*}
It is  a $1^{\rm st}$-order product form as
 $ \{0,1\} \oplus \{0, 2\} \oplus \{0, 16, 32\} \equiv {\Bbb Z}_{12}$. \qquad $\Box$

 \bigskip

 We will come back to these two examples in Section 6 to explain some situations.

\end{section}

\bigskip
\bigskip

\begin{section}{\bf $\Phi$-tree, blocking and kernel polynomials}

 In this section, we will study the mask polynomial by using a graph theoretic consideration,  and then use it to consider the condition ($P_k$). We make use of a setup by Protasov on the refinement equations [P]. Let $V_0 = \{\vartheta \}$ be the root,
$$
 V_k=\{{\bf j} = j_k\dots j_1: \ j_\ell \in \{0,1,\dots, b-1\}, \ j_1\neq 0\},  \quad k \geq 1,
 $$
 and $V=\bigcup_{k\geq 0} V_k$\ . (We reverse the usual ordering on the index as we are dealing with the integers in the $b$-adic expansion instead of the decimals.) For any  ${\bf j}\in V_k,\  k \not =0$, it has $b$ offsprings $j_{k+1}{\bf j}$  (note that by assumption, $\vartheta$ has only $b-1$ offsprings in $V$, which are the elements of $V_1$). Let $E$ be the set of edges connecting those ${\bf j}$ and $j_{k+1}{\bf j}$.  Then $(V,E)$ is a tree with $\vartheta$ as the root, and  we call it  a \emph{Protasov tree} (associated with $b$).

For each ${\bf j} \in V_k$, we let \ $
m_{\bf j}  = j_kb^{k-1} + \cdots + j_2 b + j_1$, which is  the  $b$-adic expansion determined by ${\bf j}$. Note that $j_1 \not =0$ by the assumption on the Protasov tree, it follows that there is a one-to-one correspondence between $V_k$ and the set of integers in $\{1 , \cdots ,  b^k-1\}$ which are not divisible by $b$.

\bigskip

We call $B \subset V\setminus \{\vartheta \}$ a \emph{blocking} if it is a finite set and every infinite path starting from $ \vartheta $ must intersect exactly one element of $B$. The following criterion is due to Protasov \cite{[P]} on the refinement equation adjusted to the present situation. For each ${\bf j} \in V_k$, we use $e(m_{\bf j})$ to denote  $e^{2\pi i m_{\bf j}/b^k}$ for simplicity.

\bigskip

\begin{theorem} \ \label{th7.4} ${\mathcal D}$ is a tile digit set if and only if there is a blocking $B$ such that
for any  ${\bf j}\in B$,
\begin{equation}\label{7.2}
P_{\mathcal D}\big (e (m_{\bf j})\big ) =0.
\end{equation}
(We call such $B$ a $P_{\mathcal D}$-blocking)
\end{theorem}

\bigskip

\noindent {\bf Remark}. \    The Kenyon criterion involves  checking all integers $m>0$, while the Protasov criterion only involves checking finitely many $m$'s (although the finding of the blocking set $B$ is not direct). The seemingly weaker tree criterion actually implies the Kenyon criterion in the following way: Suppose $B$ is  such a blocking. Let $m\geq 1$ be an integer such that $b\nmid m$, we write $m$ in $b$-adic expansion as
$m=m_{\bf j}$ with ${\bf j} = j_t\cdots j_1 \in V$.  If  $ {\bf j} $ has an ancestor  ${\bf i } = j_k \cdots j_1$ belonging to $B$, then we choose this $k$ for the $b^k$ in the
Kenyon criterion; otherwise  there is an $\ell$ so that  ${\bf i} = \underbrace {0\cdots0}_{\ell}  j_t\dots j_1 \in B$,  then $k=\ell+t$ satisfies the Kenyon criterion.

\bigskip

 In the following we will convert the Protasov tree into a tree of cyclotomic polynomial, which is more tractable to study the structure of the tile digit sets.

 \bigskip

 For  ${\bf j} \in V_k$, we let $a_{{\bf j}}=\hbox{g.c.d.} (m_{{\bf j}},b^k)$,  $d = d_{{\bf j}} ={b^k}/a_{{\bf j}}$, and associate with ${\bf j}$  a  $\Phi_d$  so that $\Phi_d (e^{2\pi i {m_{\bf j}}/{b^k}})=0$. In this way, we define a map $\tau$ from $V$ to the set of all cyclotomic polynomials by mapping ${\bf j}$ to $\Phi_{d_{\bf j}}$. (By convention $\tau (\vartheta) = \vartheta$.)

 Let $\Phi_d$ in the range of $\tau$, we define
 $$
  C_d =  \{{\bf i}: \tau({\bf i}) = \Phi_d\} = \tau^{-1}(\Phi_d);
 $$
 if  ${\bf j}\neq\vartheta$ is such that $\tau({\bf j}) = \Phi_d$, we define
 $$
 L_{{\bf j}} = \{\ell{\bf j}: 0\leq\ell\leq b-1\} \  \mbox{and } \ C_d^*= \bigcup_{{\bf i}\in C_d} L_{{\bf i}}.
 $$
 For  ${\bf i} = i_k...i_1$, we denote ${\bf i}^- = i_{k-1}...i_1$ . The following proposition is some basic properties of the map $\tau$.

 \bigskip

 \begin{Prop} \label{lem7.1} With the above notations, we have

 \vspace {0.2cm}

(i) If $\tau({\bf i}) = \tau({\bf j})$, then ${\bf i}$ and ${\bf j}$ lie in the same $V_k$, and
$\tau({\bf i}^-) = \tau({\bf j}^-)$;

\vspace {0.3cm}

(ii) $\#C_d = \deg\Phi_d$, and
\begin{equation}\label{7.1}
\Phi_d(x) = {\prod}_{{\bf j} \in C_d} \big(x- e({m_{\bf j}})\big);
\end{equation}

\vspace {0.3cm}
(iii) If $ \tau({\bf j}) = \Phi_d $, then $\tau(L_{{\bf j}}) = \{\Phi_e: \Phi_e(x)|\Phi_d(x^b)\}$ and
\begin{equation}\label{7.1'}
\Phi_d(x^b)={\prod}_{{\bf j}' \in C^*} \big (x-e(m_{{\bf j}'})\big ).
\end{equation}
 \end{Prop}

 \bigskip

 \noindent {\bf Proof}.
 (i) Suppose that ${\bf i} \in V_{k_1}$, ${\bf j} \in V_{k_2}$ with $k_1<k_2$ and $\tau({\bf i}) = \tau({\bf j}) =\Phi_d$. Let  $a_{\bf i} = \hbox {g.c.d.}(b^{k_1},m_{\bf i}), \ a_{\bf j} = \hbox {g.c.d.}(b^{k_2},m_{\bf j})$, then
$$
d \ = \ {b^{k_1}}/{a_{\bf i}}\  = \ {b^{k_2}}/a_{\bf j}.
$$
Hence $a_{\bf j} =b^{k_2-k_1}a_{\bf i}$,  so that $b$ divides $m_{\bf j}$, which implies $j_1 =0$.  This contradicts   the assumption that $j_1 \not =0$ for ${\bf j} \in V$.

To prove the second part, we let ${\bf i} = i_k \cdots i_1$ and ${\bf j} = j_k \cdots j_1$. Since $\tau({\bf i}) = \tau({\bf j})$, we have  $\hbox {g.c.d.}(b^k,m_{\bf i} )= \hbox {g.c.d.}(b^k, m_{\bf j})$. This is equivalent to
$$
m_{\bf j}=c_1b^k+c_2m_{\bf i} \quad \text{ and } \quad m_{\bf i}=c_1'b^k+c_2'm_{\bf j}
$$
for some integers $c_1,c_2,c_1',c_2'$. Hence
$$
m_{\bf j}=(c_1b)b^{k-1}+c_2m_{\bf i}\quad  \text{ and }\quad   m_{\bf i}=(c_1'b)b^{k-1}+c_2'm_{\bf j},
$$
so that $\hbox {g.c.d.}(b^{k-1}, m_{\bf i})=\hbox {g.c.d.}(b^{k-1},m_{\bf j})$. Consequently
$ \hbox {g.c.d.}(b^{k-1},m_{{\bf i}'}) = \\ \hbox {g.c.d.}(b^{k-1},m_{{\bf j}'}).$

\bigskip

(ii) By (i),  $C_d\subset V_k$ for some $k$. Hence $C_d = \{{\bf j}\in V_k : \Phi_{d}(e(m_{{\bf j}}))=0\}$, so that $\#C_d\leq \deg\Phi_d$. On the other hand, let $a_{\bf j} = {\hbox {g.c.d.}(b^k,m_{{\bf j}})}$,  then $d = {b^k}/ a_{\bf j}$  for some ${\bf j}\in C$. Note that $\deg\Phi_d$ equal the Euler-phi function $\varphi (d)$, i.e.,
\begin {equation} \label {eq6.1+}
\deg\Phi_d = \#\{r: 1\leq r< d, \ \hbox {g.c.d.}(r,d)=1\}.
\end {equation}
 Hence for each such $r$, we have $r< d = {b^k}/a_{\bf j}$. This implies $r a_{\bf j} < b^k$ and we can write $r a_{\bf j} = m_{{\bf j}'}$ for a unique ${\bf j}'\in V_k$. In this case, $ e^{2\pi i r/d} = e(m_{{\bf j}'})$.  We have shown  that each $r$ in \eqref{eq6.1+} corresponds to exactly one ${\bf j}' \in V_k$ with $\tau ({\bf j}') = \Phi_d$. Hence $\deg\Phi_d\leq \#C_d$.

\medskip

  For ${\bf j} \in C_d \subset V_k$, $ e(m_{\bf j}) :=e^{2\pi i m_{\bf j}/{b^{k}}}$ is a root of $\Phi_d(x)$. This implies that the product  in \eqref {7.1} is a factor of $\Phi_d(x)$. Since  $\deg\Phi_d =  \#C_d$, the polynomials involved are monic,  identity  (\ref {7.1}) must hold.

\medskip

(iii) Given ${\bf j}$ such that $\tau({\bf j}) = \Phi_d$.  by noting that
$e^{2\pi i (\ell b^k +m_{\bf j})/b^{k+1}}$ with $0\leq \ell \leq b-1$ and ${\bf j} \in C_d$ are roots  of $\Phi_d(x^b)$, we have  $\tau(L_{{\bf j}}) \subset \{\Phi_e: \Phi_e(x)|\Phi_d(x^b)\}$.

\medskip

Conversely, let $a_{\bf j} = \hbox {g.c.d.}(b^{k},m_{\bf j})$, we have $m_{{\bf j}} = ra_{{\bf j}}$, $b^{k} = da_{{\bf j}}$ and $r,d$ is relatively prime. We then write $b = sd'$ where  $d'$ contains all prime factors appearing in $d$. We note that
\begin{equation}\label{7.3}
\frac{\ell b^k+ m_{{\bf j}}}{b^{k+1}} = \frac{\ell da_{{\bf j}} + ra_{{\bf j}}}{bda_{{\bf j}}} = \frac{\ell d +r}{bd} = \frac{\ell d +r}{sd'd}
\end{equation}
 and $\hbox {g.c.d.}(\ell d +r, d'd)=1$. Now, given any $\Phi_e(x)|\Phi_d(x^b)$, $e = s'dd'$ where $s'$ is a factor of $s$ (by Proposition \ref{Prop3.1}(ii)). Since $\hbox {g.c.d.}(d,s)=1$, $\{\ell d+ r: 0\leq \ell\leq b-1\}$ contains a complete residue set (mod $s$). This implies that there exists $\ell_0$ such that $\ell_0d+r\equiv s/s' (\mbox{mod} \ s)$ and so that  $\frac{\ell_0 d+r}{s} = \frac{c}{s'}$ with $\hbox {g.c.d.}(c,s')=1$.  By \eqref {7.3},
$$
\frac{\ell_0 b^k+ m_{{\bf j}}}{b^{k+1}} = \frac{c}{s'd'd}
$$
and $\hbox {g.c.d.}(c, s'd'd)=1$. This implies that $\tau(\ell_0 m_{\bf j}) = \Phi_e$, and $\{\Phi_e: \Phi_e(x)|\Phi_d(x^b)\}\subset \tau(L_{{\bf j}})$ follows.

\medskip

The  identity \eqref{7.1'} follows from
$e^{2\pi i (\ell b^k +m_{\bf j})/b^{k+1}}$ with $0\leq \ell \leq b-1$ and ${\bf j} \in C_d$ are roots  of $\Phi_d(x^b)$, and the degree of the two polynomials in the identity is
$b \deg (\Phi_d)$.
 \qquad $\Box$

\bigskip

 By this proposition,  we can associate a Protosav tree with a {\it tree of cyclotomic polynomials (with respect to b)},  which we  call it a \textit{$\Phi$-tree} (see Figure 1). The set of vertices of this tree at level $1$ are $\tau(V_1)$.  The offsprings of $\Phi_d$ are the cyclotomic  factors of $\Phi_d(x^b)$, they are determined by Proposition \ref {Prop3.1}(ii). By (iii), an edge joining ${\bf j}$ to its offspring $\ell{\bf j}$ corresponds to an edge joining $\Phi_d$ and a cyclotomic  factor of $\Phi_d(x^b)$. We can see easily from (iii) that the map $\tau$ is surjective on the $\Phi$-tree.  Moreover, all $\Phi_d$ in the tree are different by Proposition \ref{Prop3.3}.

\bigskip

\begin{figure}[h]
\centerline{\includegraphics[width=12cm,height=4.5cm]{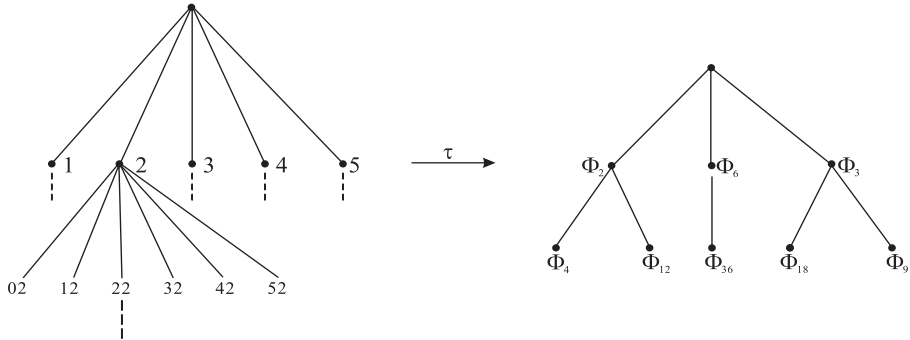}}
\caption {\small{ An illustration of the Protasov tree and the associated $\Phi$-tree  for the case $b=6$. On the first level,  $\tau (1)= \tau(5)= \Phi_6;  \ \tau(2) = \tau(4) = \Phi_3; \  \tau(3) = \Phi_2$. On the second level, $\Phi_3(x^6) = \Phi_{9}(x)\Phi_{18}(x)$, hence $\Phi_3$ has two descendants $\Phi_{9}$ and $\Phi_{18}$. The map $\tau$ acts on the descendant of $2$ as:  $\tau(02) =\tau(22) =\tau(32) = \Phi_9; \ \tau(12) = \tau(32) = \tau (52) = \Phi_{18}$.}}
\end{figure}

\bigskip

A blocking $B$ in $V$  is called  a \emph{symmetric blocking } if  $\tau^{-1}(\tau({\bf j}))\subset B$ for every  ${\bf j}\in B$.

 \bigskip

 \begin{Lem}\label{lem7.3+}
 Let  $B$ be a symmetric  blocking in the Protasov tree, then ${\mathcal N}_B: = \tau (B)$ is a blocking in the $\Phi$-tree. Conversely, if ${\mathcal N}$ is a blocking in the $\Phi$-tree, then $B_{{\mathcal N}} = \tau^{-1}({\mathcal N})$ is a symmetric blocking in the Protasov tree.
 \end{Lem}

 \noindent{\bf Proof.} It is clear that ${\mathcal N}_B$ is a finite set. Let $\Gamma = \{\vartheta, \Phi_{d_1},\Phi_{d_2},\cdots\}$ be an infinite path in the $\Phi$-tree (i.e. $\Phi_{d_k}(x)$ is a factor of $\Phi_{d_{k-1}}(x^b)$). Since $B$ is symmetric, by Proposition \ref{lem7.1}(iii), there exists an infinite path $ \gamma = \{\vartheta, {\bf j}_1,{\bf j}_2,\cdots \}$ in the Protasov tree such that $\tau({\bf j}_k) = \Phi_{d_k}$ for all $k$.  Since $B$ is a blocking, we can find a vertex ${\bf j}_r$ on $\gamma$ such that ${\bf j}_r\in B$. Hence $\Phi_{d_r}\in {\mathcal{N}}_B$. That $\Gamma$ hits ${\mathcal{N}}_B$ exactly once is an easy consequence from the blocking $B$ in $V$.

 \bigskip

Conversely, let $\gamma  =\{\vartheta, {\bf j}_1,{\bf j}_2,...\}$ be an infinite path in $V$, then Proposition \ref{lem7.1}(iii) implies that $\tau(\gamma) = \{\vartheta ,\tau({\bf j}_1),...\}$ is a path in the $\Phi$-tree. Hence, there is a unique $\tau({\bf j}_k) = \Phi_{d}\in {\mathcal N}\cap \tau(\gamma)$. It follows that  ${\bf j}_k\in B_{{\mathcal N}}$, i.e., $\gamma$ meets $B_{{\mathcal N}}$ at ${\bf j}_k$ and it is the unique point by Proposition \ref{lem7.1}(i). This shows that  $B_{{\mathcal N}}$ is a blocking, and it is clearly  symmetric.
  \qquad $\Box$

 \medskip

 \begin{Def}\label{Def7.6}  Let ${\mathcal N}$ be a blocking of the $\Phi$-tree, we call
$$
K(x) =\prod_{\Phi_d\in{\mathcal{N}}} \Phi_d(x)
$$
a kernel polynomial (with respect to $b$).
\end{Def}

\bigskip

Note that the most basic kernel polynomial is
$
K(x) = \prod_{d|b, d>1} \Phi_d(x).
$
We can generate all the kernel polynomials by adopting the following simple procedure step by step,  starting from $ K_{ {\tau (\mathcal E)}}(x)$: let $K_{\mathcal N}$ be as in the definition, and let $\Phi_d \in {\mathcal N}$. The descendants of $\Phi_d$, denoted by ${\mathcal N}_d$,  are the factors of $\Phi_d(x^b)$. It follows that ${\mathcal N}' = ({\mathcal N} \setminus \{\Phi_d\}) \cup {\mathcal N}_d $ is again a blocking, and the corresponding kernel polynomial is
\begin{equation} \label {eq7.3}
K_{{\mathcal N}'}(x) = \frac{\Phi_d(x^b)K_{\mathcal N}(x)}{ \Phi_d(x)} \ .
\end{equation}

\bigskip

\begin{Lem}\label{lem7.5}  A minimal $P_{\mathcal D}$-blocking (i.e., smallest cardinality) is symmetric.
\end{Lem}

\medskip

\noindent {\bf Proof}. Let $B$ be a minimal $P_{\mathcal D}$-blocking and let ${\bf j}$ be an element of $B$.
Suppose there exists ${\bf i} $ such that $\tau({\bf i}) = \tau({\bf j})$ but ${\bf i}\not\in B$. If there is no ancestor of ${\bf i}$ belongs to $B$, then we can form a $B'$ by including  ${\bf i}$ into $B$ and drop all its offsprings in $B$.  This  $B'$  is a $P_{\mathcal D}$-blocking since $\tau({\bf i}) = \tau({\bf j})$,  and has smaller cardinality than $B$, it is impossible by the minimality assumption. Hence there is an ancestor  ${\bf i}'$ of ${\bf i}$  such that ${\bf i}'\in B$. Let ${\bf j}'$ be the ancestor of ${\bf j}$ such that $|{\bf j}' |=|{\bf i}'|$, then by Proposition \ref{lem7.1}(i), $\tau({\bf i}') = \tau({\bf j}')$. By dropping all the offsprings of ${\bf j}'$ and adding ${\bf j}'$ to B, we obtain a  smaller $P_{{\mathcal D}}$-blocking. We see that it is again impossible, and the lemma follows.
\qquad $\Box$

\bigskip

Our main conclusion is the following characterization of the tile digit sets.
\medskip

\begin{theorem}\label{th7.9} Let $b>1$ be an integer and let ${\mathcal{D}}$ be a digit set with $\#{\mathcal D}=b$. Then the following are equivalent.

\ \ (i) ${\mathcal D}$ is a tile digit set of $b$;

\ (ii) there is a symmetric $P_{\mathcal D}$-blocking in the Protasov tree of $b$;

(iii) there is  a blocking ${\mathcal N}$ in the $\Phi$-tree of $b$ such that
$$
K(x) : = \prod_{\Phi_d\in{\mathcal N}}\Phi_d(x)
$$
\hspace{1.0cm} is a kernel polynomial  and  $ K(x) | P_{\mathcal{D}}(x)$;

(iv) $P_{\mathcal{D}}(x)$ satisfies condition $(P_k)$  for some $k\geq 1$.

\end{theorem}

\medskip

\noindent {\bf Proof}.
$(i)\Rightarrow(ii)$. Since we can always find a minimal blocking, the assertion follows from the Protasov criterion (Theorem \ref{th7.4}) and Lemma \ref{lem7.5}.

\medskip
$(ii)\Rightarrow(iii)$. For a symmetric blocking $B$ in (ii), Lemma \ref{lem7.3+} implies that ${\mathcal{N}}_B$ is a blocking in the $\Phi$-tree. For any $\Phi_d\in{\mathcal N}_{B}$, $\Phi_d = \tau({\bf j})$ for some ${\bf j}\in B$. By (ii), $P_{\mathcal D}\big (e(m_{\bf j})\big ) =0$, this shows that $\Phi_d(x)|P_{\mathcal{D}}(x)$.  Since ${\mathcal{N}}_B$ is a blocking, $K(x)$ is a kernel polynomial by definition, and  clearly it divides $P_{\mathcal{D}}(x)$.

\medskip

 $(iii)\Rightarrow(iv)$.  We can localize the $(P_k)$ condition by checking $(P_{k_d})$ for each $d|b$, then take $k = \max \{k_d : \  d|n\}$.  Note that by assumption, $P_{\mathcal D} (x) = K(x)Q(x)$, hence it suffices to prove that for each $d|b$, $K(x)$ satisfies some $(P_{k_d})$ condition.

 We trace through the steps of producing $K(x) := K_{\mathcal N}(x)$ as in \eqref{eq7.3}. For a fixed $d|b$, there exists a smallest $j_1$ such that some factors of $\Phi_d(x^{b^{j_1}})$ divide $K(x)$ (as ${\mathcal N}$ is a blocking). If all factors of $\Phi_d(x^{b^{j_1}})$ divides $K(x)$, then $K(x)$ satisfies  $(P_1)$ for $d$.

Otherwise, there exists at least one factor in $\Phi_{d}(x^{b^{j_1}})$ that does not divide $K(x)$, say $e_1$.  We repeat the same procedure on $\Phi_{e_1}$ as for $\Phi_d$ and find $j_2$ such that some factors in $\Phi_{e_1}(x^{b^{j_2}})$ will divide $K(x)$. We check whether all factors of $\Phi_{e_1}(x^{b^{j_2}})$ will divide $K(x)$, i.e,  to satisfy condition $(P_2)$. If not, we continue on with the same procedure. Since the blocking ${\mathcal N}$ is a finite set, the process must stop finally. This factor $d$ must satisfy some $(P_{k_d})$ condition.

\medskip

$(iv)\Rightarrow(i)$. This has been proved in Proposition \ref{Prop6.1}.
\qquad $\Box$

\bigskip

\noindent {\bf Remarks}. \

(1)\ For the $1^{\rm st}$-order product-form ${\mathcal D}$, we have seen that $K_{{\mathcal D}}$ defined in (\ref{5.7-}) is a polynomial satisfying $(P_1)$ and it divides $P_{{\mathcal D}}(x)$, it is  a kernel polynomial for ${\mathcal D}$. Similarly, for the ${\mathcal D}$ in \eqref{eq5.3}, the polynomial $K_{{\mathcal D}}$ is also a kernel polynomial of ${\mathcal D}$, it satisfies $(P_2)$.

\medskip


\medskip
In Section 4, we have shown that a $k^{\rm th}$-order modulo product-form satisfies the $(P_k)$ condition. Actually for $b=p^{\alpha}q$, the converse is also true. This is contained in our classification of the kernel polynomials for $ b= p^\alpha q$ in [LLR].

\bigskip

(2)\   For a tile digit set ${\mathcal D}$, $P_{\mathcal D}$ may contain more than one kernel polynomials. For example, let $P_{\mathcal D}(x)$ be the mask polynomial of ${\mathcal D}$ in Example \ref{ex6.3}.  It contains the factors $\Phi_{2^73^3}(x)$, $\Phi_{2^83^3}(x)$ and $\Phi_{2^6\cdot3^2}(x)$. We write
$$
P_{\mathcal D} (x) = K_1(x) \Phi_{2^73^3}(x) \Phi_{2^6\cdot3^2}(x) = K_2(x)\Phi_{2^73^3}(x)\Phi_{2^83^3}(x)
$$
with $K_1(x), K_2(x)$ defined in an obvious way. It is seen that $K_1(x) \Phi_{2^73^3}(x) $ is the mask polynomial of the tile digit set in Example \ref{ex6.2}.  This $K_1(x)$ forms a blocking in the $\Phi$-tree of $b=12$  and is thus a kernel polynomial of ${\mathcal D}$ (the blocking property can be checked easily using the diagram of Figure 2 below). Moreover, note that both factors $\Phi_{2^6\cdot3^2}(x)$ and $\Phi_{2^83^3}(x)$ are on the same path in the $\Phi$-tree, and the path has no branches.  We can interchange the two factors to form $K_2(x)$ which is still a blocking, and therefore is another kernel polynomial of ${\mathcal D}$.

\medskip

\begin{figure}[h]
\centerline{\includegraphics[width=8cm,height=3cm]{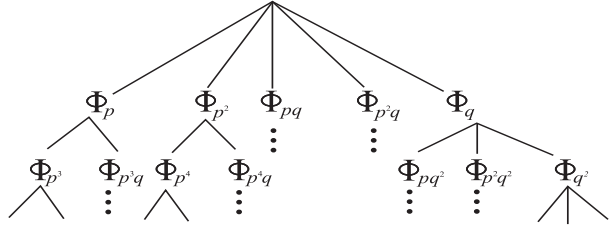}}
\caption {{ \small The $\Phi$-tree for  $p^2q$. }}
\end{figure}

\bigskip

(3) \   It is possible that a kernel polynomial $K(x)$ does not generate any tile digit set at all. For example, let $b=pq$, consider
\begin{equation}\label{6.3}
 K(x) = \Phi_{p}(x)\Phi_{q}(x)\Phi_{p^kq^k}(x), \quad  k \geq 1.
\end{equation}
It is clear that $K(x)$ is a blocking (see Figure 1). For $k=1$, $K(x)$ corresponds to $\{0, \cdots , b-1\}$.  However for $k > 1$, $K(x)$ cannot generate a tile digit set (i.e., $K(x)$ does not divide $P_{\mathcal D}(x)$ for any tile digit set ${\mathcal D}$.   To see this, we note that from [LR] the weak product-form characterization implies that the kernel polynomial can only be the form $\Phi_p(x)\Phi_{q^{\ell}}(x^{p^{\ell}})$. Now from (\ref{6.3}) the prime power spectrum is $\{p,q\}$. Hence, $\ell=1$ and  $k$ in (\ref{6.3}) can only be 1 in order to obtain a tile digit set for $b=pq$.

\medskip

  It is not clear which kernel polynomials can generate a tile digit set. A complete answer is known only for $b=p^{\alpha}$ (every kernel polynomial can generate tile digit sets) and $p^{\alpha}q$, and a partial answer is also given when $b=p^{\alpha}q^{\beta}$. These details are given in  \cite{[LLR]}.  In here we provide a simple sufficient condition to this question.

\bigskip

\begin{Prop}\label{prop7.8} Let $K(x)$ be a kernel polynomial associated with $b$, let $Q(x)$ be a polynomial in ${\mathbb Z} [x]$ with $Q(1)=1$. If $P(x)=K(x)Q(x)$ has non-negative coefficients, then $P(x) = P_{\mathcal D}(x)$ for some tile digit set ${\mathcal D}$.
\end{Prop}

\noindent {\bf Proof}.  Write $P(x)$ explicitly as
$$
P(x)=a_0+ a_1 x^{d_1}+\dots +a_m x^{d_m},
$$
where $a_j\geq 1$ are integers and $\sum_{\ell=1}^m a_\ell=P(1)=K(1)Q(1)=b$. Hence $m\leq b$.

Let ${\mathcal D}=\{0=d_0,d_1,\dots,d_m\}$. We claim that $m=b$. Suppose $m<b$, we note that $(b, {\mathcal D})$ define an IFS (as in (\ref {eq1.0})) and a self-similar set $T(b, {\mathcal D})$.
Let
$$
p_j=  \frac {a_j}{{\sum}_{\ell=1}^m a_\ell}
$$
be a set of probability on the IFS. It induces a self-similar measure $\mu$ with  supported on  $T(b, {\mathcal D})$ [H]. Note that $P(x)$ is the mask polynomial of $\mu$, and $K(x)$ is a kernel polynomial. They imply that  $P(x)$ satisfies the Protasov criterion for a refinement equation [P].  It follows that  $\mu$ is  absolutely continuous, and its support $T(b, {\mathcal D})$  has positive Lebesgue measure.  However,  $\#{\mathcal D}=m<b$ implies that $T(b, {\mathcal D})$ has zero Lebesgue measure. This is a contradiction. Hence $m=b$. This means $a_j$ must all be $1$ and ${\mathcal D}$ must be a tile digit set with $P(x) = P_{\mathcal D}(x)$.
\qquad $\Box$

\bigskip

\end{section}

\bigskip

\bigskip

\begin{section} {\bf Remarks and open questions}

Our consideration in this paper is in terms of tile digit sets, in fact, many of the theorems also hold for the absolutely continuous self-similar measures. Let ${\mathcal D}\subset {\Bbb Z}^+$ such that ${\mathcal D} = m $ (may not equal $b$), and  let $\phi_{d}(x) = b^{-1}(x+d)$ be the corresponding affine map, it is well known that there exists a unique probability measure of equal weight $\mu = \mu(b,m,{\mathcal D})$ which satisfies
\begin{equation}\label{6.2}
\mu (E) = \sum_{d\in{\mathcal D}} \frac{1}{m}\mu (\phi_d^{-1}(E)),
\end{equation}
for any Borel set $E$.

\medskip

A natural question is to determine when $\mu$ is absolutely continuous with respect to the Lebesgue measure. It is known that the general criterion of Kenyon and Protasov for the mask polynomial of ${\mathcal D}$ is necessary and sufficient for $\mu$ to be absolutely continuous ([DFW], [P]). We can modify Theorem \ref{th4.2} and Theorem \ref{th7.9} as follows:

\medskip

{\it Let $\mu = \mu(b,{\mathcal D})$ be the self-similar measure defined in (\ref{6.2}). Then

 (i) If $\mu$ is absolutely continuous with respect to the Lebesgue measure, then Theorem \ref{th4.2} holds except for the uniqueness of prime-power spectrum,

\vspace {0.2cm}

 (ii) If $\mu$ is absolutely continuous, then the equivalence of Theorem \ref{th7.9} holds. }

\bigskip

It is easy to show that if the above $\mu$ is absolutely continuous, then $m \geq b$. As an application, we prove the following interesting proposition.

\medskip

 \begin{Prop}
 If $\mu = \mu(b,m,{\mathcal D})$ is
absolutely continuous, then  $m = bk$ for some integer $k\geq 1$.
\end{Prop}

\medskip
\noindent {\bf Proof}.
Let $b=p_1^{\alpha_1}...p_k^{\alpha_k}$ be its unique prime factorization.
If $\mu(b,m,{\mathcal D})$ is absolutely continuous, then Theorem \ref{th4.2} shows that for each $p_j$ there exists $\{a_{j,\ell}\}_{\ell = 0}^{\alpha_j-1}$ with $a_{j,\ell}\equiv\ell (\mbox{mod} \ \alpha_j )$ and $\Phi_{p_j^{a_{j,\ell}}}(x)$  divides $P_{\mathcal D}(x)$. Let
$$
S(x) = \prod_{j=1}^{k}\prod_{\ell=0}^{\alpha_{j}-1}\Phi_{p_j^{a_{j,\ell}}}(x).
$$
We have $P_{\mathcal D}(x) = S(x)Q(x)$ for some integral polynomial $Q(x)$. Note that the cyclotomic factors are all distinct, by Proposition \ref{Prop3.1}(iii) we have  $S(1) = p_1^{\alpha_1}...p_k^{\alpha_k} =b$. Hence,
$$
 m= P_{\mathcal D}(1) = S(1)Q(1)
$$
and $Q(1)$ is an integer since $Q(x)$ is integral polynomial. This completes the proof.
\qquad $\Box$

\bigskip

   There are also many unanswered basic questions on the self-affine tiles and the structure of the tile digit sets.  So far most of our consideration is on ${\Bbb R}^1$, and very little is known on ${\Bbb R}^s$. One of the difficulties is that there is no analogue for the cyclotomic polynomials in higher dimension. In [HL], it was proved that for $|\det (A)| = \#{\mathcal D} =p$ a prime, if the linear span  span$({\mathcal D})$ is ${\Bbb R}^s$, then $T(A, {\mathcal D})$ is a self-affine tile if and only if ${\mathcal D}$ is a complete residue set with respect to $A$.

 \bigskip

 \noindent\textbf{Q1:}  Can the condition span$({\mathcal D})={\Bbb R}^s$ be omitted in the above characterization?

 \bigskip

  It is easy to generalize the product-form to higher dimension [LW3].  It will be interesting to set up the following:

  \bigskip

 \noindent\textbf{Q2:} Find the corresponding generalization for condition ($P_k$) and the modulo product-forms in ${\Bbb Z}^s$, and  use them to study the tile digit sets.

\bigskip

Our study on the  structure  of the tile digit sets is closely related to the spectral set problems. Recall that a closed subset $\Omega \subset {\Bbb R}^s$ is called a {\it spectral set} if $L^2(\Omega)$ admits an exponential orthonormal basis $\{e^{2\pi i \langle \lambda, \cdot \rangle}\}_{\lambda\in\Lambda}$. The well-known Fuglede conjecture asserted that $\Omega$ is a spectral set if and only if it is a translational tile. Recently Tao \cite{[T]} gave a counter-example that the conjecture fails in ${\Bbb R}^s$ with $s\geq 5$, and  the example was eventually modified to $s\geq 3$ [KM]. Now the conjecture remains open in ${\Bbb R}^1$ and ${\Bbb R}^2$ (see e.g.,  [\L], [{\L}aW]). For the more restrictive class of self-affine tiles, the following is still open:

\bigskip

 \noindent\textbf{Q3:} Suppose  $T(A, {\mathcal D})$ is a self-affine tile in ${\Bbb R}^s$, is $T(A, {\mathcal D})$  a spectral set?

\bigskip
\bigskip

\noindent{\bf Acknowledgement.}  The authors would like thank Professor Xinggang He for suggesting the reference [P]. They also  thank Professor Dejun Feng for bringing their attention to the absolute continuity of the self-similar measures.
\end{section}

\end{document}